\documentclass{article}
\usepackage{amsmath}
\usepackage{amssymb}
\usepackage{amsopn}

\begin{document}
\title{Gap Solitons in Almost Periodic One-Dimensional Structures}
\author{Alexander Pankov\\
Department of Mathematics\\ Morgan State University\\ 1700 E. Cold Spring Lane\\ Baltimore, MD 21251, USA}
\date{}

\maketitle

\begin{abstract}
We consider almost periodic stationary  nonlinear Schr\"odinger equations in dimension $1$. Under certain assumptions
we prove the existence of nontrivial finite energy solutions  in the strongly indefinite  case. The proof
is based on a carefull analysis of the energy functional restricted to the so-called generalized Nehari manifold, and
the existence and fine properties of special Palais-Smale sequences. As an application, we show that certain one
dimensional almost periodic photonic crystals possess gap solitons for all prohibited frequencies.

\vspace{1ex}


\vspace{1ex}

\end{abstract}

\newtheorem{theorem}{Theorem}[section]
\newtheorem{lemma}{Lemma}[section]
\newtheorem{definition}{Definition}[section]
\newtheorem{example}{Example}[section]
\newtheorem{remark}{Remark}[section]
\newtheorem{prop}{Proposition}[section]
\newtheorem{cor}{Corollary}[section]

\numberwithin{equation}{section}

\section{Introduction}\label{s1}

In this paper we consider the problem of existence of non-zero finite energy solutions (also known as bound states, or
homoclinics) to the following
one-dimensional stationary nonlinear Schr\"odinger equation (NLS)
$$
-u''+V(x)u=f(x,u)
$$
in which the $x$-dependence is almost periodic, while the linear part of the equation is not nesessarily positive
definite. More precisely, we suppose that $0$ is not in the spectrum of the linear part. The most interesting case is
when $0$ belongs to a finite spectral gap, {\em i.e.}, there is a non-empty  part of the spectrum below zero. It is
well-known that the spectrum of a  periodic Schr\"odinger operator is absolutely continuous and  has the so-called
band-gap structure. Moreover, typical one-dimensional Schr\"odinger operators have infinitely many gaps \cite{re-si-3}.
In the almost periodic case the spectrum is not absolutely continuous in general. However, it possesses gaps. Moreover,
typically the spectrum of an almost periodic Schr\"odinger operator is nowhere dense (see, {\em e.g.}, \cite{b-t,
pas-fi}).

\vspace{1ex}

In last decades, the periodic NLS  in arbitrary dimension  has been studied extensively including strictly indefinite
case (see, {\em e.g.}, \cite{di, kr-sz, pa-umzh, pa-mil, ra, szu-we, wi} and references therein).  In the almost
periodic case the situation is totally different. The first result in this direction obtained by variational methods
concerns second order Hamiltonian systems with positive definite linear part \cite{se-ta-te}, including
one-dimensional NLS equation. This result has been
developed in several directions, but still for problems with positive definite linear part (see \cite{al-cal, cz-m-n,
m-n-t, ra-ap-1, ra-ap-2, sprad, te}).

One of the key ingredients in \cite{se-ta-te} is the construction of special Palais-Smale
sequences, known as $(\overline{PS})$ sequences \cite{c-z-e-s}, based on mountain pass geometry and (negative) gradient
flow of the associated functional $J$. Notice that
the mountain pass minimax class is invariant with respect to standard deformations and, hence, the gradient flow. In
the strictly indefinite case the
 functional $J$ possesses infinite dimentional linking geometry \cite{kr-sz, wi}. However, the minimax class related to
this geometry is not invariant with respect to the gradient flow. We overcome this difficulty by employing the
generalized Nehari manifold of the functional $J$ in its original version introduced in \cite{pa-mil}. Special
Palais-Smale sequences are then constructed directly via the negative gradient flow of the functional $J$ restricted to
the generalized Nehari manifold. Notice that this requires certain additional smoothness of the nonlinearity with
respect to $u$ to guarantee the existence and uniquness for such flow.

An essential part of \cite{se-ta-te} is devoted to detailed structure of Palais-Smale sequences with the aim to relate
special Palais-Smale sequences and returning sequences of real numbers for the functional $J$. The arguments are quite
involved and depend crucially on the positivity of the linear part. In our work we restrict ourselves to Palais-Smale
sequences at levels close to the ground level, which is the infimum of $J$ over the Nehari manifold. The structure of
such sequences is not complicated so that to pass to a returnig sequence it is enough to use relatively
simple concentration-compactness arguments.

In addition, let us point out that in this paper we use a weaker concept of almost periodicity, the so-called Stepanov
almost periodicity. We do that to allow piece-wise continuous dependence of the potential $V(x)$ and nonlinearity
$f(x,u)$ on $x$. This is important in the application of our result to nonlinear optics.

\vspace{1ex}

Now let us turn to applications. The term {\em gap soliton} was born in the area of photonic crystals.
Photonic crystals are optical media with spatially periodic, or close-to-periodic, structure. Here close-to-periodic can
be almost periodic, or asymptotically periodic, or something simillar. In this context almost periodicity models
disordered periodic structures, while asymptotic periodicity represents a localized defect in a periodic structure. One
of the basic fiture of photonic crystals is
that light of certain frequences (so-called prohibited frequences) can not propagate through such a medium. This is due
to the band-gap structure of the spectrum of, say, periodic Maxwell operators. Actually, prohibeted frequences are
exactly the poins in gaps of spectrum. However, if a photonic crystal is made of non-linear media, a completely new
phenomenon occur. In such crystals there may exist localized light pattern with prohibited carier frequences. These are
called gap solitons. For physics and mathematics of photonic crystals we refer to \cite{ak-jo, joan, ku, mills, pa-mil,
npc} and references therein.

Gap solitons are widely studied in physics literature by means of numerical and asymptotical methods. However, to the
best of our knowledge there is only one regorous mathematical result on the existence of gap solitons \cite{pa-mil}. It
concerns gap solitons of special form (the so-called $TM$-mode) in one- and two-dimensional periodic crystals. As we
will see below, our main result provides the existence of gap solitons in one-dimensional almost periodic photonic
crystals.

Notice that the stationary NLS also apears as an equation for the profile function of a standing wave in the
evolutionaly nonlinear Schr\"odinger equation. Typically, such standing waves exist if its frequency belongs to a
spectral gap of the linear part. Often such waves are also called gap solitons. Certainly, our result provides the
existence of such waves under appropriate assumptions.

\vspace{1ex}

The organization of the paper is as follows. Section~\ref{s2} contains certain facts on one-dimentional Schr\"odinger
operators and reminds the concept of Stepanov almost periodicity. In Section~\ref{s3} we formulate  our main result,
while Section~\ref{s3a} is devoted to a variational formulation of the problem and certain simple results on the
continuous dependence of the energy functional on the envelope of the problem. Sections~\ref{s4} and \ref{s5} form a
core of our techniques, and are devoted to the generalized Nehari manifold and Palais-Smale sequences, respectively.
The proof of main result is contained in Section~\ref{s6}. Finally, in Section~\ref{s7}, we sketch an application to
photonic crystals.

\section{Preliminaries}\label{s2}

First, let us introduce basic spaces of real valued functions on $\mathbb{R}$.

\smallskip

 By $L^2(\mathbb{R})$ we denote the space of square integrable functions endowed with the standard norm $\|\cdot\|_2$
and
inner product $(\cdot,\cdot)$. The Sobolev space
$$
H^1(\mathbb{R})=\{u\in L^2(\mathbb{R})\, |\, u'\in
L^2(\mathbb{R})\}
$$
with the graph norm $\|\cdot\|$ is a Hilbert space. The inner product in $H^1(\mathbb{R})$ is denoted
by $\langle\cdot,\cdot\rangle$. By $L^\infty(\mathbb{R})$ we denote the space of all essentially bounded functions with
the standard norm $\|\cdot\|_\infty$.  The space of all infinitely differentiable compactly supported functions is
denoted by $C^\infty_0(\mathbb{R})$.

 By $H^{-1}(\mathbb{R})$ we denote the dual space to $H^1(\mathbb{R})$ with the norm $\|\cdot\|_*$. The symbol $(\cdot,
\cdot)$ stands  both for the inner product
in $L^2(\mathbb{R})$ and for the duality pairing on $H^{-1}(\mathbb{R})\times H^1(\mathbb{R})$. This does not lead to
any  confusion. It is well-known that
$$
H^1(\mathbb{R})\subset L^2(\mathbb{R})\subset H^{-1}(\mathbb{R})
$$
continuously and densely. Moreover, $H^1(\mathbb{R})$ is continuously embedded into $L^\infty(\mathbb{R})$. Actually,
any $H^1$-function is continuous and vanishes at infinity.

\smallskip

A locally integrable function $u$ is {\em Stepanov bounded} if
$$
\|u\|_{BS}=\sup_{t\in\mathbb{R}}\int_t^{t+1}|u(x)|\,dx <\infty\,.
$$
Such functions form a Banach space denoted by $BS(\mathbb{R})$. A function $u\in BS(\mathbb{R})$ is {\em Stepanov almost
 periodic} if the set of its shifts
$$
\{T_zu\}_{z\in\mathbb{R}}\,,
$$
where $(T_z)u(x)=u(x+z)$, is precompact in the space
$BS(\mathbb{R})$. In other words, for any sequence $z_k\in\mathbb{R}$ there exists a subsequence $z_{k'}$ such that the
sequence $T_{z'_k}u$ converges in the space $BS(\mathbb{R})$. The space of Stepanov almost periodic functions is a
closed subspace of $BS(\mathbb{R})$ denoted by $S(\mathbb{R})$. For a Stepanov  almost periodic function $u$,  the
closure of $\{T_zu\}_{z\in\mathbb{R}}$ in  the space $BS(\mathbb{R})$ is denoted by $\mathcal{E}(u)$ and is called the
{\em envelop} of $u$. The following simple fact is well-known (see, {\em e.g.}, \cite{le-zhi, pa-ap}). If $u_h=\lim
T_{z_k}u\in \mathcal{E}(u)$, then $u=\lim T_{-z_k}u_h$ (limits in the space $BS(\mathbb{R})$). The set $\mathcal{E}(u)$
 is a compact set in $BS(\mathbb{R})$. Notice that the operators $T_z$ form a strongly continuous group of operators in
 $S(\mathbb{R})$, but this is not so in the whole space $BS(\mathbb{R})$.

\smallskip

Let $V\in BS(\mathbb{R})$. Then the operator
\begin{equation}\label{e2.1}
L=L_0+V(x)=-\frac{d^2}{dx^2}+V(x)\,,
\end{equation}
defined by means of the sum of quadratic forms associated to $L_0$ and $V$,  is a bounded below self-adjoint operator in
$L^2(\mathbb{R})$. The form domain of $L$ is the space $H^1(\mathbb{R})$.
Furthermore, the operator $L$ extends to a bounded linear operator from $H^1(\mathbf{R})$ into $H^{-1}(\mathbb{R})$
still  denoted by $L$, and the extension depends continuously on $V\in BS(\mathbb{R})$ with respect to the operator
norm, hence, with respect to the norm resolvent convergence (see, {\em e.g.}, \cite{si}). Furthermore, the operator of
multiplication by $V$ is a bounded linear operator from $H^1(\mathbf{R})$ into $H^{-1}(\mathbb{R})$ and its norm does
not exceed $\|V\|_{BS}$, {\em i.e.},
\begin{equation}\label{e2.2}
|(Vu,v))|=|\int_{\mathbb{R}}V(x)u(x)v(x) dx|\leq \|V\|_{BS}\|u\|\|v\|\,.
\end{equation}
Moreover, this operator is $L_0$-form bounded with form bound $0$ \cite{si}.

\vspace{1ex}

In what follows we denote by $\sigma(L)$ the spectrum of  $L$. If $0\not\in\sigma(L)$, we denote by
$E^+\subset H^1(\mathbb{R})$ and $E^-\subset H^1(\mathbb{R})$ the positive and negative subspaces of the form $(Lu,u)$,
respectively. These subspaces are orthogonal with respect to both $H^1$ and $L^2$ inner products. Moreover, $L E^\pm$ is
orthogonal to $E^\mp$ with respect to duality pairing $(\cdot,\cdot)$ on $H^{-1}(\mathbb{R})\times H^1(\mathbb{R})$. By
$P^+$ and $P^-$ we denote the orthogonal projectors in $H^1(\mathbb{R})$ onto $E^+$ and $E^-$, respectively. Notice that
these projectors are orthogonal with respect to $(\cdot,\cdot)$ as well. Each element $u\in H^1$ possesses the
representation $u=u^++u^-$, where $u^+=P^+u$ and $u^-=P^-u$.

\begin{prop}\label{p2.1}
Let $V\in BS(\mathbb{R})$. If $0\not\in\sigma(L)$, then there exists a constant $\kappa>0$, depending on $\|V\|_{BS}$
and
the distance between $0$ and $\sigma(L)$, such that
\begin{equation}\label{e2.4}
(Lu,u)\geq\kappa\|u\|^2\,, \quad u\in E^+\,,
\end{equation}
and
\begin{equation}\label{e2.5}
(Lu,u)\leq -\kappa\|u\|^2\,, \quad u\in E^-\,,
\end{equation}
\end{prop}

{\em Proof\/}. We prove inequality (\ref{e2.4}), the other is similar.

Let $2\delta$ be the distance between zero and $\sigma(L)$. Then
$$
(Lu,u)\geq 2\delta\|u\|^2_{2}\,,\quad u\in E^+\,.
$$
Since $V$ is $L_0$-form bounded with form bound $0$, then there exist sufficiently small $\alpha\in (0, 1)$ and
$\beta>0$,
depending on $\|V\|_{BS}$, such that
$$
|(Vu,u)|\leq \alpha (L_0u,u)+\beta\|u\|^2_{2}=\alpha\|\dot u\|^2_{2}+\beta\|u\|^2_{2}\,,\quad u\in H^1(\mathbb{R})\,.
$$
Hence, for all $u\in E^+$,
$$
(1-\alpha)\|u\|^2\leq (Lu,u)+ C\|u\|^2_{2}\,,
$$
where $C=1+\beta-\alpha$. The right hand side of this inequality can be expressed as
$$
\frac{C+\delta}{\delta}\left[\frac{\delta}{C+\delta}((Lu,u)-\delta\|u\|^2_{2})+\delta\|u\|^2_{2}\right]\,.
$$
Since $\delta/(C+\delta)<1$, on the subspace $E^+$ this quantity does not exceed
$$
\frac{C+\delta}{\delta}\,(Lu,u)\,,
$$
and the result follows.

\hfill$\Box$

Now suppose that $V\in S(\mathbb{R})$.  The {\em envelop} $\mathcal{E}(L)$ of $L$ is the set of all operators $L_h$ of
the form (\ref{e2.1}) generated by potentials $V_h\in \mathcal{E}(V)$. Being considered as a subset in the Banach space
of all bounded linear operators from $H^1(\mathbb{R})$ into $H^{-1}(\mathbb{R})$, the envelop $\mathcal{E}(L)$ is a
compact set.

\begin{prop}\label{p2.2}
Suppose that $V\in S(\mathbb{R})$. Then $\sigma(L_h)=\sigma(L)$ for all $L_h\in\mathcal{E}(L)$.
\end{prop}

{\em Proof\/}. If $V_h\in\mathcal{E}(V)$, then there exists a sequence $z_k\in\mathbb{R}$ such that $T_{z_k}V\to
V_h$ in $BS(\mathbb{R})$. It is easily seen that $\sigma(L_0+T_{z_k}V)=\sigma(L)$. Since $L_0+T_{z_k}V\to
L_h$ with respect to the norm resolvent convergence, then, by \cite[Theorem~VIII.23]{re-si-1},
$\sigma(L)\subset\sigma(L_h)$. But $T_{-z_k}V_h\to V$ in $BS(\mathbb{R})$. Hence, interchanging the role of
$V$ and $V_h$, we obtain the required.

\hfill$\Box$

\begin{remark}\label{r2.1}
If $V\in S(\mathbb{R})$ and $0\not\in\sigma(L)=\sigma(L_h)$, we denote by $E^+_h$ and $E^-_h$ the positive and negative
subspaces of the quadratic form $(L_hu,u)$. By Proposition~\ref{p2.2}, the conclusion of Proposition~\ref{p2.1} holds
for $L_h$ with the same constant $\kappa$. Furthermore, positive and negative spectral projectors depend continuously on
the potential. More precisely, let $T_{z_k}V\to V_h$ in $BS(\mathbb{R})$, and let $P^\pm_k$ and
$P^\pm_h$ be the positive (negative) spectral projector that correspond to the potentials $T_{z_k}V$ and
$V_h$, respectively. Then $P^\pm_k\to P^\pm_h$ with respect to the operator norm.
\end{remark}

For functions of two variables, $g(x,u)$, we need an appropriate concept of almost periodicity with respect to the first variable $x\in\mathbb{R}$. It is always assumed that such a function is a Carath\'eodory function, {\em i.e.}, $g(x,u)$ is continuous in $u$ for almost all $x\in\mathbb{R}$, and  Lebesgue measurable in $x$ for all $u\in\mathbb{R}$. For any $R>0$, we set
$$
\|g\|_R=\|\sup_{|u|\leq R}|g(\cdot, u)|\|_{BS}\,.
$$
We say that $g(x,u)$ is {\em strictly Stepanov almost periodic} in $x$ (in symbols $g\in S(\mathbb{R}\times\mathbb{R})$) if $\|g\|_R<\infty$ for all $R>0$, and for any sequence $z_k\in\mathbb{R}$ there exist a subsequence $z_{k'}$ and a function $g_h$ such that $\|g_h\|_R<\infty$ for all $R>0$ and
$$
\|T_{z'_k}g(\cdot, u)-g_h(\cdot, u)\|_R\to 0 \quad \forall R>0\,.
$$
In other words, being considered as a function of $x\in\mathbb{R}$ with values in the (Frech\'et) space of continuous functions of $u\in\mathbb{R}$, $g$ is a Stepanov almost periodic function. The envelope $\mathcal{E}(g)$ of $g$ consists of all such limit functions $g_h$. Notice, that any strictly Stepanov almost periodic function is Stepanov almost periodic in $x$ uniformly with respect to $u\in [-R,R] \,\,\, \forall R>0$, but not vise versa.

\section{Statement of Problem and Main Result}\label{s3}

We are looking for nonzero vanishing at infinity  solutions to the following one-dimensional nonlinear Schr\"odinger
equation
\begin{equation}\label{e3.1}
- u''(x)+V(x)u(x)=\chi f(x,u(x))\,,
\end{equation}
where $\chi=\pm 1$.

Let
$$
F(x,u)=\int_0^u f(x,s)\,ds\,.
$$
Throughout  the remaining part of the paper we suppose that the following assumptions hold true.
\begin{description}
\item[$(i)$] {\em The potential $V$ is Stepanov almost periodic, $V\in S$,  and the spectrum of the operator $L$ does
not
contain zero. In the case when $\chi=-1$ we suppose in addition that there is a non-empty part of the
spectrum below $0$\/}.

\item[$(ii)$] {\em For almost all $x\in\mathbb{R}$, the function $f(x,u)$ is  continuously differentiable with respect to $u\in\mathbb{R}$.  The functions $F(x,u)$, $f(x,u)$ and $f_u(x,u)$ are strictly Stepanov almost periodic. For any $u\neq 0$,  the function $F(x,u)$ is bounded below by a positive constant\/}.

\item[$(iii)$] {\em The nonlinearity satisfies  $f(\cdot,0)=0$ and $f_u(\cdot,0)=0$. Furthermore, for every $R>0$ there  
exists a constant $\mu(R)>0$ such that
$$
|f_u(x,u)-f_u(x,v)|\leq \mu(R)|u-v|\,,\quad  |u|\,,|v|\leq R\,.
$$
for almost all $x\in\mathbb{R}$}.

\item[$(iv)$] {\em There exists a constant $\theta\in (0,1)$ such that for almost all $x\in\mathbb{R}$}
$$
0<f(x,u)u\leq \theta\cdot f_u(x,u)u^2\,,\quad  u\neq 0\,.
$$
\end{description}

Without loss of generality we suppose that  $\mu(R_1)\leq \mu(R_2)$  whenever $R_1\leq R_2$.

Assumption $(i)$ guarantees that the self-adjoint operator $L$ is well-defined (see Section~\ref{s2}). By the mean value theorem, Assumption $(iii)$ implies that for almost all $x\in\mathbb{R}$
\begin{equation}\label{e3.2-1}
|f(x,u)|\leq  \mu(R)|u|^2
\end{equation}
and
\begin{equation}\label{e3.3-1}
|F(x,u)|\leq \mu(R)|u|^3
\end{equation}
whenever $|u|\leq R$.
Assumption $(iv)$ implies easily that
\begin{equation}\label{e3.3a}
0< qF(x,u)\leq f(x,u) u\,,\quad  u\neq 0\,,
\end{equation}
where $q=(1+\theta)/\theta>2$. This is the standard Ambrosetti-Rabinowitz condition. In particular, from (\ref{e3.3a}) it follows  that for any $\varepsilon>0$ there exists a constant $C_\varepsilon>0$ such that
\begin{equation}\label{e3.3b}
F(x,u)\geq -\varepsilon |u|^2+C_\varepsilon|u|^q\,.
\end{equation}

Notice that in Assumption $(ii)$ it is enough to assume strict Stepanov almost periodicity for $f_u$ only. Then so is for $f$ and $F$.

\vspace{1ex}

{\bf Example.} The nonlinearity
\begin{equation}\label{e-2-ex}
f(x,u)=\alpha(x)|u|^{p-2}u\,,
\end{equation}
satisfies Assumptions $(ii)$--$(iv)$ provided $\alpha\in S(\mathbb{R})\cap L^\infty(\mathbb{R})$, $\mathrm{ess} \inf\alpha>0$, and $p\geq 3$.

\vspace{1ex}

Under Assumptions imposed above,  the set of shifts $\{(T_zV, T_zf)\}_{z\in \mathbb{R}}$ is precompact with repect to
the topology generated by semi-norms $\|V\|_{BS}+\|f_u\|_R$, $R>0$. Its closure is denoted by $\mathcal{E}$. This is a
compact set. In what follows we always suppose that the set $\mathcal{E}$ is parameterized, not necessarily in a
one-to-one way,
by elements $h$ of an index set $\mathcal{H}\supset \mathbb{R}$. Together with equation (\ref{e3.1}) we consider the
following family of  equations
\begin{equation}\label{e3.4}
-u''(x)+V_h(x)u(x)=\chi f_h(x,u(x))\,,\quad h\in\mathcal{E}\,.
\end{equation}
These equations form the {\em envelop} of equation (\ref{e3.1}), which can be identified with $\mathcal{E}$.
Any equation in the envelop satisfies Assumptions $(i)$--$(iv)$ with the same $\mu(R)$ and $\theta$.

\vspace{1ex}

Our main result is the following.

\begin{theorem}\label{t3.1}
Under Assumptions $(i)$--$(v)$ equation (\ref{e3.1}) has a nonzero solution $u\in H^1(\mathbb{R})$. Moreover, the
solution $u$ is continuously differentiable and decays at infinity exponentially fast, {\em i.e.}, there exist positive
constants $\alpha$ and $\beta$ such that
$$
|u(x)|+|u'(x)|\leq\alpha\exp(-\beta|x|)\,.
$$
\end{theorem}

The solution in Theorem~\ref{t3.1} is a weak solution, {\em i.e.\/},
$$
\int_{\mathbb{R}}(u'(x)\varphi'(x)+V(x)u(x)\varphi(x))\,dx=\chi\int_{\mathbb{R}}f(x,u(x))\varphi(x)\,dx
$$
for all $\varphi\in C^\infty_0(\mathbb{R})$.

\begin{remark}\label{r3.1}
Theorem~\ref{t3.1} applies to all equations (\ref{e3.4}) in the envelop of equation (\ref{e3.1}).
\end{remark}

\begin{remark}
Suppose that  zero is below the essential spectrum of $L$, {\em i.e.,} $L$ is positive definite. If $\chi=-1$, then
it is easily seen that equation (\ref{e3.1}) has only trivial solution in $H^1(\mathbb{R})$. If $\chi=1$ and $V(x)\geq
\alpha_0>0$, the existence of nontrivial solution is obtained for a wider class of nonlinearities, including
(\ref{e-2-ex}) with $p>2$ (see \cite{se-ta-te}). Actually, in \cite{se-ta-te} the potential is a constant function,
while $f(x,u)$
is Bohr almost periodic in $x$, but the arguments of that paper extend straightforwardly to the case of non-constant
potential and Stepanov almost
periodic $x$-dependence.
\end{remark}

\section{Variational Formulation}\label{s3a}

Associated to equation (\ref{e3.1}), we introduce the functional
\begin{equation}\label{3.5}
\begin{split}
J(u)&=\frac{1}{2}\int_{\mathbb{R}}(|u'(x)|^2+V(x)u^2(x)\,dx)-\chi\int_{\mathbb{R}}F(x,u(x))\,dx\\
&=\frac{1}{2}(Lu,u)-\chi \Phi(u)\,.
\end{split}
\end{equation}

Similarly, we introduce the functional $J_h$ associated to equation (\ref{e3.4}). Its non-quadratic part is denoted by $\Phi_h$. The functionals $J_h$ form the envelop of $J$.

Under the assumptions imposed above,
the functional $J$ is a well-defined $C^{2,1}$-functional on the space $H^1(\mathbb{R})$. Its first and second derivatives are given by
\begin{equation}\label{e3.6}
(J'(u),v)=(Lu,v)-\chi\int_{\mathbb{R}} f(x,u(x))v(x)\,dx\,,\quad u,v\in H^1(\mathbb{R})\,,
\end{equation}
and
\begin{equation}\label{e3.7}
(J''(u)v,w)=(Lv,w)-\chi\int_{\mathbb{R}}f_u(x,u(x))v(x)w(x)\,dx\,,\quad u,v,w\in H^1(\mathbb{R})\,.
\end{equation}
Notice that $J'$ is weakly continuous.

Often it is convenient to use gradients of $J$ instead of derivatives. These are defined by
$$
\langle\nabla J(u),v\rangle=(J'(u),v)
$$
and
$$
\langle\nabla^2 J(u)v,w\rangle=(J''(u)v,w)
$$
for all $u, v, w\in H^1(\mathbb{R})$. Then $\nabla J(u)\in H^1(\mathbb{R})$, while $\nabla^2J(u)$ is a linear bounded operator in $H^1(\mathbb{R})$.

\vspace{1ex}

Now we estimate the difference between two functionals of the form $J_h$ and its derivative.

\begin{prop}\label{p3.1}
For any $h_i\in\mathcal{H}$, $i=1, 2$, and any $R>0$
$$
|J_{h_1}(u)-J_{h_2}(u)|\leq \frac{1}{2}\|V_{h_1}-V_{h_2}\|_{BS}\|u\|^2+\|(f_{h_1})_u-(f_{h_2})_u\|_R\|u\|^2
$$
and
$$
\|J'_{h_1}(u)-J'_{h_2}(u)\|_*\leq \|V_{h_1}-V_{h_2}\|_{BS}\|u\|+\|(f_{h_1})_u-(f_{h_2})_u\|_R\|u\|
$$
provided $u\in H^1(\mathbb{R})$ with $\|u\|\leq R$.
\end{prop}

{\em Proof\/}. By inequality (\ref{e2.2}), both the difference of the linear parts and its derivative are estimated by the first term in the right hand sides.

By Taylor's formula and inequality (\ref{e2.2})
\begin{equation*}\begin{split}
|\Phi_{h_1}(u)-\Phi_{h_2}(u)|&\leq \int_{\mathbb{R}}\int_0^1|(f_{h_1})_u(x, tu(x))-(f_{h_1})_u(x, tu(x))|(1-t)u^2(x) dt dx\\
&\leq\int_{\mathbb{R}}\sup_{|u|\leq R}|(f_{h_1})_u(x, u)-(f_{h_1})_u(x, u)|u^2(x) dx\\
&\leq \|(f_{h_1})_u-(f_{h_1})_u\|_R\|u\|^2
\end{split}\end{equation*}
which implies the first estimate of the proposition.

The proof of second inequality is similar.

\hfill$\Box$

\begin{prop}\label{p3a.1}
If $u_n\to u_0$ weakly in $H^1(\mathbb{R})$, then
\begin{equation}\label{e3a.1}
J_h(u_n-u_0)-J_h(u_n)+J_h(u_0)\to 0
\end{equation}
and
\begin{equation}\label{e3a.2}
J'_h(u_n-u_0)-J'_h(u_n)+J'_h(u_0)\to 0
\end{equation}
strongly in $H^{-1}(\mathbb{R})$ uniformly with respect to $h\in\mathcal{H}$.
\end{prop}

{\em Proof\/}. The integrand of non-quadratic part, $\Psi_h$, of $J_h$
satisfies inequalities (\ref{e3.2-1}) and (\ref{e3.3-1}) uniformly with respect
to $h\in\mathcal{H}$. Hence, arguing exactly as in \cite[Lemma 2.6]{te} we
obtain the result of proposition for $\Psi_h$ instead of $J_h$.
Due to linearity of the operator $L_h$, this implies (\ref{e3a.2}) immediately.

The quadratic part of the left-hand side in (\ref{e3a.1}) coincides with
$$
(L_hu_0,u_0)-(L_hu_0,u_n)\to 0\,,
$$
and we obtain (\ref{e3a.1}) for every individual $h\in\mathcal{H}$. Since the operators $L_h$ form a compact set of
bounded linear operators from $H^1(\mathbb{R})$ into $H^{-1}(\mathbb{R})$, this convergence is uniform with respect to
$h\in\mathcal{H}$.

\hfill$\Box$

The proof of Theorem~\ref{t3.1} is given in the subsequent sections. Obviously,
$u=0$ is a trivial critical point of the functional  $J$.  We shall prove that
$J$ possesses a nontrivial critical point. {\em In the course of the proof we
consider equation (\ref{e3.1}) in the case when $\chi=1$\/}.
The other case is completely similar. We only need to replace the functional $J$
by $-J$ and interchange the role of the subspaces $E^+$ and $E^-$ introduced in
Section~\ref{s2}.

\vspace{1ex}

\section{Generalized Nehari Manifold}\label{s4}

The {\em generalized Nehari manifold} $\mathcal{N}$ of the functional $J$ consists of all nonzero $u\in H^1(\mathbb{R})$ such that
$$
(J'(u),u)=0
$$
and
$$
(J'(u),v)= 0\,, \quad \forall v\in E^-\,.
$$
Equivalently, these equations can be written as $\langle\nabla J(u),u\rangle=0$ and  $P_-\nabla J(u)=0$, respectively.
 The generalized Nehari manifold of a functional $J_h\in\mathcal{E}$ is denoted by $\mathcal{N}_h$.

For any $w\not\in E^-$ we set
$$
E_w=\{sw+v : s > 0, v\in E^-\}
$$
and
$$
\bar{E}_w=\{sw+v : s\in\mathbb{R}, v\in E^-\}\,.
$$
By the definition of $\mathcal{N}$, if $u$ is a critical point of $J|_{E_w}$, then $u\in\mathcal{N}$. As consequence, $\mathcal{N}$ contains all nontrivial critical points of $J$.

\begin{lemma}\label{l4.1} For every  $w\not\in E^-$, the functional $J|_{E_w}$
attains its positive global maximum.
\end{lemma}

{\em Proof\/}. Without loss of generality, we can suppose that $w\in E^+$ and $\|w\|=1$. If $s\in (0,1]$, then, by (\ref{e3.3-1}),
$$
J(sw)\geq\frac{s^2}{2}(Lw,w)-\mu_1s^3\,.
$$
Hence, $J(sw)>0$ for $s>0$ small enough.

On the other hand, by (\ref{e2.5}) and (\ref{e3.3b}), for any $sw+v\in E_w$
$$
J(sw+v)\leq\frac{1}{2}s^2(Lw,w)-\frac{1}{2}\kappa\|v\|^2+\varepsilon s^2\|w\|^2_{L^2}+\varepsilon\|v\|^2_{L^2}-C_\varepsilon\|sw+v\|^q_{L^q}\,.
$$
Since the norm of a projector in a Banach space is $\geq 1$, we have that 
$$
\|sw+v\|_{L^q}\geq C\|sw\|_{L^q}\,.
$$ 
Then
$$
J(sw+v)\leq(\frac{1}{2}(Lw,w)+\varepsilon\|w\|^2_{L^2})s^2 -(\frac{1}{2}\kappa-\varepsilon)\|v\|^2-C'_\varepsilon\|w\|_{L^q}^qs^q\,.
$$
Taking $\varepsilon$ small enough, we obtain that $J(sw+v)\to-\infty$ as $\|sw+v\|\to\infty$.

Obviously, $J|_{E_w}$ is upper weakly semi-continuous. Hence, it attains its (positive) global maximum.

\hfill$\Box$

\begin{remark}\label{r4.0} As in \cite[Proposition 2.3]{szu-we}), one can show that
for every $w\not\in E^-$ the intersection $\mathcal{N}\cap E_w$ consists of exactly one point which is a  unique 
maximum point of $J|_{E_w}$. But we do not use this fact.
\end{remark}

It is convenient to introduce the functional
$$
I(u)=J(u)-\frac{1}{2}(J'(u),u)\,.
$$
Obviously, $J(u)=I(u)$ for all $u\in \mathcal{N}$. By inequality (\ref{e3.3a}), $I(u)\geq 0$ for all $u\in H^1(\mathbb{R})$.

Now we prove the following technical result.

\begin{lemma}\label{l4.2}
There exists a constant $C>0$ independent of $u\in H^1(\mathbb{R})$ such that
\begin{equation}\label{e4.1}
\|u\|^2\leq C(|(J'(u),u)|+|(J'(u),u^-)|+\mu(\|u\|_\infty)\|u\|_\infty\|u\|^2)
\end{equation}
and
\begin{equation}\label{e4.2}
\|u\|^2\leq C(|(J'(u),u)|+|(J'(u),u^-)|+(I^{1/2}(u)+I(u))\|u\|)
\end{equation}
for all $u\in H^1(\mathbb{R})$.

\end{lemma}

{\em Proof\/}. The identity
$$
(J'(u),u^-)=(Lu^-,u^-)-\int_\mathbb{R} f(x,u)u^- dx
$$
and Proposition~\ref{p2.1} imply
\begin{equation}\label{e4.3}
\kappa\|u^-\|^2\leq -(J'(u),u^-)-\int_\mathbb{R}f(x,u)u^- dx\,.
\end{equation}
Similarly,  the identity
$$
(J'(u),u)=(Lu^+,u^+)-\int_\mathbb{R} f(x,u)u^+ dx +(J'(u),u^-)
$$
implies
\begin{equation}\label{e4.4}
\kappa\|u^+\|^2\leq (J'(u),u) -(J'(u),u^-)+\int_\mathbb{R}f(x,u)u^+ dx\,.
\end{equation}
Adding inequalities (\ref{e4.3}) and (\ref{e4.4}), we obtain immediately that
\begin{equation}\label{e4.5}
\begin{split}
\|u\|^2&\leq C (|(J'(u),u)|+|(J'(u),u^-)|+\\
&+\int_\mathbb{R}|f(x,u)||u^+| dx+\int_\mathbb{R}|f(x,u)||u^-| dx)\,.
\end{split}
\end{equation}
Then, by inequality (\ref{e3.2-1}),
$$
\int_\mathbb{R}|f(x,u)||u^\pm| dx\leq\mu(\|u\|_\infty)\|u\|_\infty\int_{\mathbb{R}}|u||u^\pm| dx\leq \mu(\|u\|_\infty)\|u\|_\infty\|u\|_2\|u^\pm\|_2\,.
$$
Hence,
\begin{equation*}
\begin{split}
\|u\|^2&\leq C (|(J'(u),u)|+|(J'(u),u^-)|+\mu(\|u\|_\infty)\|u\|_\infty\|u\|_2(\|u^+\|_2+\|u^-\|_2))\\
&\leq C (|(J'(u),u)|+|(J'(u),u^-)|+\mu(\|u\|_\infty)\|u\|_\infty\|u\|^2)\,,
\end{split}
\end{equation*}
which proves (\ref{e4.1}).

\vspace{1ex}

Now we prove inequality (\ref{e4.2}). Given $u\in H^1(\mathbb{R})$, let
$$
S_1=\{x\in\mathbb{R} : |u(x)|\leq 1\}
$$
and $S_2=\mathbb{R}\setminus S_1$. We introduce the following integrals
$$
I_1=\int_{S_1}|f(x,u)|^2 dx
$$
and
$$
I_2=\int_{S_2}|f(x,u)| dx\,.
$$

By inequality (\ref{e3.2-1}), $f^2(x,u)\leq\mu(1)f(x,u)u$ on $S_1$, while on $S_2$ we have  that  $|f(x,u)|\leq
f(x,u)u$.
Then, by inequality (\ref{e3.3a}),
\begin{equation}\label{e4.6}
I(u)\geq (2^{-1}-q^{-1})\int_{\mathbb{R}}f(x,u)u dx\geq\nu I_k\,, \quad k=1, 2\,,
\end{equation}
for some $\nu>0$.
Since
\begin{equation*}
\begin{split}
\int_\mathbb{R}|f(x,u)||u^\pm| dx &\leq  (\int_{S_1}|f(x, u)|^2 dx)^{1/2}(\int_{S_1}|u^\pm|^2 dx)^{1/2} +\\
&+\|u^\pm\|_\infty\int_{S_2}|f(x, u)| dx\\
&\leq (I_1^{1/2}+I_2)\|u^\pm\|\,,
\end{split}
\end{equation*}
equations (\ref{e4.5}) and (\ref{e4.6}) yield (\ref{e4.2}).

\hfill$\Box$

\begin{prop}\label{p4.1}
There exists a constant $\varepsilon_0>0$ such that $\|u\|\geq\|u\|_\infty\geq\varepsilon_0$, $J(u)\geq\varepsilon_0$ and
$$
\int_\mathbb{R} f(x,u)u dx\geq 2\varepsilon_0
$$
 for all $u\in\mathcal{N}$.
\end{prop}

{\em Proof\/}. The first two statements follow immediately from Lemma~\ref{l4.2}. Since $F\geq 0$, we see that $L(u)\geq 2\varepsilon_0$ on $\mathcal{N}$. Now the last statement follows from the definition of $\mathcal{N}$.

\hfill$\Box$

\begin{remark}\label{r4.1}
Obviously, Lemma~\ref{l4.2} and Proposition~\ref{p4.1} hold for all functionals $J_h$ in the envelop of $J$ with the same constants $C$ and $\varepsilon_0$. In particular, for any nontrivial critical point $u$ of $J_h$ we have that $\|u\|\geq\|u\|_\infty\geq\varepsilon_0$ and $J_h(u)\geq\varepsilon_0$.
\end{remark}

Let ${\bar E}^-=\mathbb{R}\oplus E^-$. Elements of this space are denoted by $[\tau, v]$, where $\tau\in\mathbb{R}$ and $v\in E^-$.  The inner product in this space is still denoted by $\langle\cdot,\cdot\rangle$. We introduce the operator $G : H^1(\mathbb{R})\rightarrow {\bar E}^-$ by the formula
$$
G(u)=[\langle \nabla J(u), u\rangle,\, P^-\nabla J(u)]\,,\quad u\in H^1(\mathbb{R})\,.
$$
It is not difficult to verify that the operator $G$ is a $C^{1,1}$ map, and its derivative is given by the formula
$$
G'(u)v=[\langle\nabla^2 J(u)v,u\rangle+\langle\nabla J(u), v\rangle,\, P^-\nabla^2 J(u)v]
$$
for all $u, v\in H^1(\mathbb{R})$. Notice that $\mathcal{N}=G^{-1}(0)\setminus\{0\}$.

\begin{lemma}\label{l4.3}
Let $u_0\in H^1(\mathbb{R})$, and let
$$
\gamma_0=\langle\nabla J(u_0),u_0\rangle
$$
and
$$
\gamma=P^-\nabla J(u_0)\in E^-\,.
$$
 Then, for all $\tau\in\mathbb{R}$ and $v\in E^-$,
\begin{equation}\label{e4.7}\begin{split}
\langle G'(u_0)(\tau u_0+v),[\tau, v]\rangle&\leq 2\gamma_0\tau^2-\kappa\|v\|^2+\frac{3}{2}\tau^2\|\gamma\|\\
&+\frac{3}{2}\|\gamma\|\|v\|^2-\tau^2(1-\theta)\int_{\mathbb{R}}f(x, u_0)u_0 dx\,,
\end{split}
\end{equation}
where $\kappa>0$ and $\theta\in (0, 1)$ are constants from Proposition~\ref{p2.1} and Assumption $(iv)$, respectively.
\end{lemma}

{\em Proof\/}. Since, by the assumptions,
$$
(Lu_0,u_0)=\gamma_0+\int_{\mathbb{R}}f(x,u_0)u_0 dx
$$
and
$$
(Lu_0,v)=\langle\gamma, v\rangle +\int_{\mathbb{R}} f(x, u_0)v dx\,,
$$
a straightforward, but a little bit tedious, calculation yields the identity
\begin{equation}\begin{split}\label{e4.8}
\langle G'(u_0)(\tau u_0+v),[\tau, v]\rangle&=2\tau^2\gamma_0+(Lv,v)+3\tau\langle\gamma, v\rangle\\
&- \int_{\mathbb{R}}(H(x)\tau^2+2K(x)\tau v+M(x)v^2) dx\,,
\end{split}
\end{equation}
where
$$
H(x)=f_u(x,u_0)u^2_0-f(x,u_0)u_0\,,
$$
$$
K(x)=f_u(x,u_0)-f(x,u_0)
$$
and
$$
M(x)=f_u(x,u_0)\,.
$$
Obviously,
$$
|\tau\langle\gamma, v\rangle|\leq \frac{1}{2}\|\gamma\|(\tau^2+\|v\|^2)
$$
and, by Proposition~\ref{p2.1},
$$
(Lv,v)\leq -\kappa\|v\|^2\,.
$$
Therefore, it is enough to show that
$$
H(x)\tau^2+2K(x)\tau v(x)+M(x)v^2(x)\geq \tau^2(1-\theta)f(x,u_0(x))u_0(x)\,.
$$
Notice that this inequality is trivial for all $x\in\mathbb{R}$ such that $u_0(x)=0$. Suppose now that $u_0(x)\neq 0$. In this case $M(x)\neq 0$, and
\begin{equation*}\begin{split}
H\tau^2+2K\tau v+Mh^2&=\left(H-\frac{K^2}{M}\right)\tau^2+\left(\sqrt{M}v+\frac{K}{\sqrt{M}}\right)^2\\
&\geq \left(H-\frac{K^2}{M}\right)\tau^2\,.
\end{split}\end{equation*}
Simplifying  and making use of the inequality
$$
f_u(x,u_0)\geq \theta^{-1}f_0(x, u_0)u_0^{-1}
$$
which follows from Assumption $(iv)$, we obtain that
\begin{equation*}\begin{split}
\left(H-\frac{K^2}{M}\right)&=f(x, u_0)u_0-\frac{f^2(x, u_0)}{f_u(x, u_0)}\\
&\geq (1-\theta)f(x, u_0)u_0\,.
\end{split}\end{equation*}
This implies the required.

\hfill$\Box$

\begin{lemma}\label{l4.3a}
Let $R$ be any positive number. Then

(a) For any $u_0\in\mathcal{N}$ such that $\|u_0\|\leq R$, the operator
$$
G'(u_0)|_{{\bar E}_{u_0}} : {\bar E}_{u_0}\rightarrow {\bar E}^-
$$
is invertible and the norm of its inverse operator
$[G'(u_0)|_{{\bar E}_{u_0}}]^{-1}$ is bounded above by a constant that depends on $R$ only.

(b) The norms of projectors generated by the splitting
$$
{\mathrm{ker}\, G'(u_0)} + {\bar E}_{u_0}\,,\quad u_0\in\mathcal{N}\,, \|u_0\|\leq R\,,
$$
are bounded above by a constant that depends on $R$ only.

(c) The norm
$$
\|u\|_{u_0}=\|u_1\|+\|u_2\|\,,\quad  u_0\in \mathcal{N}\,,
$$
where $u_1\in {\mathrm{ker}\, G'(u_0)}$ and $u_2\in {\bar E}_{u_0}$, is equivalent to the standard $H^1$-norm uniformly with respect to $u_0\in\mathcal{N}$ with $\|u_0\|\leq R$.
\end{lemma}

{\em Proof\/}. $(a)$ By Lemma~\ref{l4.3}, with $\gamma_0=0$ and $\gamma=0$, the composition of the isomorphism $[\tau, v]\rightarrow \tau u_0+v$ and $G'(u_0)$ is a negative definite, hence, invertible operator in ${\bar E}^-$. The norm of the inverse of above mentioned isomorphism is bounded above by a constant that depends on $R$ only.   This implies the required.

$(b)$ The projector onto ${\bar E}_{u_0}$ is given by $[G'(u_0)|_{{\bar E}_{u_0}}]^{-1}\circ G'(u_0)$. Since the operator $G'(u_0)$ is uniformly bounded while $\|u_0\|\leq R$, the result follows.

$(c)$ This is an immediate consequence of $(b)$.

\hfill$\Box$

Inspecting the standard proofs of the Inverse Function and Implicit Function theorems (see, {\em e.g.}, \cite{dr-mi}, Theorems~4.1.1 and 4.2.1), we see that the following  complements to those theorems hold true.

\begin{prop}\label{p4.1a}
 Let $\varphi : X\rightarrow Y$ be a $C^{1,1}$-map between Banach spaces such that the derivative $\varphi'$ is bounded and globally Lipschitz continuous.

$(a)$ Given $c_0>0$, there exist $\rho>0$ and $C>0$ with the following property. For every $x_0\in X$ such that $\varphi (x_0)$ is invertible and $\|\varphi'(x_0)^{-1}\|\leq c_0$, the inverse function $\varphi^{-1}$ is defined on the $\rho$-neighborhood of $\varphi(x_0)$ and its Lipschitz constant does not exceed $C$.

$(b)$ Given $c_0>0$ and $c_1>0$, there exist $\rho>0$ and $C>0$ with the following property. Let $X=X_1+X_2$ be any splitting of $X$
with mutually complementary closed subspaces as components such that
$$
\|x_1\|+\|x_2\|\leq c_1\|x_1+x_2\|
$$
for all $x_i\in X_i$, $i=1, 2$. If $x_0=x_{1,0}+x_{2,0}$  is such that $\varphi(x_0)=0$ and the partial derivative $\varphi'_2(x_0)$ along $X_2$ satisfies
$$
\|\varphi'_2(x_0)^{-1}\|\leq c_0\,,
$$
then there exists a unique $C^{1,1}$ function $\psi$ defined on the $\rho$-neighborhood of $x_{1,0}$ in $X_1$ such that
$$
\varphi(x_1+\psi(x_1))=0\,, \quad \psi(x_{1,0})=x_{2,0}\,,
$$
and the Lipschitz constant of $\psi'$ is bounded above by $C$.

\end{prop}

\begin{prop}\label{p4.2}
The set $\mathcal{N}$ is a non-empty closed $C^{1,1}$-sub-manifold of $H^1$ with the tangent space $T_{u_0}={\mathrm{ker}\, G'(u_0)}$ at $u_0\in\mathcal{N}$. Furthermore, given $R>0$, there exist $\rho>0$  and $C>0$ such that for every $u_0\in\mathcal{N}$, with $\|u_0\|\leq R$, there exists a $C^{1,1}$-diffeomorphism from the $\rho$-neighborhood of $0$ in $T_{u_0}$ onto a neighborhood of $u_0$ in $\mathcal{N}$ such that the Lipschitz constant of its derivative does not exceed $C$.
\end{prop}

{\em Proof\/}. The result follows immediately from Proposition~\ref{p4.1a}$(b)$ and Lemma~\ref{l4.3a}.

\hfill$\Box$

\begin{remark}\label{r4.3}
Since $\mathcal{N}$ is a $C^{1,1}$-manifold, its tangent spaces form a fiber bundle of class $C^{0,1}$.
\end{remark}

\begin{lemma}\label{l4.4}
Given $c_0>0$ and $c_1>0$, there exist positive numbers $\alpha=\alpha(c_0, c_1)$, $r=r(c_0,c_1)$ and $C=C(c_0,c_1)$ such that, for any $u_0\in H^1(\mathbb{R})\setminus E^+$ satisfying
$$
\int_{\mathbb{R}}f(x, u_0)u_0 dx\geq c_0\,,
$$
$\|u_0\|\leq c_1$ and $\|G(u_0)\|\leq\alpha$, the restriction
$$
G_{u_0}=G|_{\bar{E}_{u_0}} : \bar{E}_{u_0}\rightarrow \bar E^-
$$
has a local inverse $G^{-1}_{u_0}$ defined on the open ball $B(G(u_0), r)$ of radius $r$ centered at $G(u_0)$, and $\|(G^{-1}_{u_0})'(\xi)\|\leq C$ for all $\xi\in B(G(u_0), r)$.
\end{lemma}

{\em Proof\/}. By Lemma~\ref{l4.3}, given $c_0>0$ there exists sufficiently small $\alpha>0$ such that
$$
\int_{\mathbb{R}}f(x, u_0)u_0 dx\geq c_0\,,
$$
and $\|G(u_0)\|\leq\alpha$ imply that the operator $(G_{u_0})'(u_0) : \bar{E}_{u_0}\rightarrow \bar E^-$ is invertible and the norm of the inverse operator $[(G_{u_0})'(u_0)]^{-1}$ is bounded by a constant, say, $c_2>0$  that depends on $c_0$, $c_1$ and $\alpha$, hence, on $c_0$ and $c_1$ only. As consequence, there exists a local inverse map $G^{-1}_{u_0}$ in a neighborhood of $\xi_0=G(u_0)$.

Now the result follows from Proposition~\ref{p4.1a}$(a)$.

\hfill$\Box$

Let us introduce the following quantities
$$
m=\inf\{J(u)\,:\, u\in \mathcal{N}\}
$$
and
$$
m_h=\inf\{J_h(u)\,:\, u\in \mathcal{N}_h\}\,\quad h\in\mathcal{H}.
$$
By Proposition~\ref{p4.1} and Remark~\ref{r4.1}, these numbers are strictly positive.

\begin{prop}\label{p4.3}
For all functionals in the envelop of $J$ we have that $m_h=m$.
\end{prop}

{\em Proof\/}. If $h\in \mathcal{H}$, then there exists a sequence $h_k\in\mathbb{R}$ such that
$$
V_{h_k}=V(\cdot+h_k)
$$
converges to $V_h$, while
$$
F_{h_k}=F(\cdot+h_k,\cdot)\,,
$$
$$
f_{h_k}=f(\cdot+h_k,\cdot)
$$
and
$$
(f_{h_k})_u=f_u(\cdot+h_k,\cdot)
$$
converge to $F_h$, $f_h$ and $(f_h)_u$, respectively, in the sense described in Section~\ref{s2}.

Let $\varepsilon>0$, and let $u\in\mathcal{N}_h$ be such that
$$
J_h(u)\leq m_h+\varepsilon\,.
$$
Setting $u_k=u(\cdot-h_k)$, it is easily seen that $\|u_k\|=\|u\|$. By Proposition~\ref{p3.1},
$$
\int_\mathbb{R}f(x,u_k)u_k dx=\int_\mathbb{R}f(x+h_k,u)u dx\to\int_\mathbb{R}f_h(x,u)u dx\,,
$$
and
$$
J(u_k)=J_{h_k}(u)\to J_h(u)\,.
$$
In addition, making use  of the fact that the spectral projectors  depend continuously  on $h\in\mathcal{H}$ (see Remark~\ref{r2.1}), we obtain that
$$
G(u_k)\to G_h(u)=0\,,
$$
where $G_h$ is the defining operator of the manifold $\mathcal{N}_h$.

By Lemma~\ref{l4.4}, $0$ is in the domain of $G_{u_k}^{-1}$ provided $k$  is large enough. Setting ${\tilde u}_k=G_{u_k}^{-1}(0)$, we have that ${\tilde u}_k\in\mathcal{N}$ and $\|u_k-{\tilde u}_k\|\to 0$. As consequence, $J({\tilde u}_k)\to J_h(u)$. This implies immediately that $m\leq m_h+\varepsilon$. Since $\varepsilon$ is an  arbitrary positive number, $m\leq m_h$.

Interchanging the role of $J$ and $J_h$ in the previous argument, we see that $m_h\leq m$, and the proof is complete.

\hfill$\Box$

\section{Palais-Smale Sequences}\label{s5}

Remind that a {\em Palais-Smale sequence}  for the functional $J$ at level $c$ is a sequence $u_n\in H^1(\mathbb{R})$ such that $J(u_n)\to c$ and $J'(u_n)\to 0$ strongly in $H^{-1}(\mathbb{R})$ (equivalently, $\nabla J(u_n)\to 0$ strongly in $H^1(\mathbb{R})$). Also we consider Palais-Smale sequences for the restriction $J|_{\mathcal{N}}$ of the functional $J$ to the generalized Nehari manifold $\mathcal{N}$. These are defined similarly. Namely, a sequence $u_n\in\mathcal{N}$ is a Palais-Smale sequence for $J|_{\mathcal{N}}$ at level $c$ if $J(u_n)\to c$ and $\nabla_\tau J(u_n)\to 0$ strongly in $H^1(\mathbb{R})$, where $\nabla_\tau$ stands for the tangent component of the gradient.

\begin{prop}\label{p5.1}
Let $u_n\in H^1(\mathbb{R})$ be a Palais-Smale sequence for $J$ at level $c$. Then the sequence $u_n$ is bounded in $H^1(\mathbb{R})$. Furthermore, $u_n\to 0$ strongly in $H^1(\mathbb{R})$ if and only if $c=0$.
\end{prop}

{\em Proof\/}. Since $J(u_n)$ is bounded and $\|J'(u_n)\|_*\to 0$, we have that
$$
I(u_n)\leq C+\varepsilon_n\|u_n\|\,,
$$
where $\varepsilon_n\to 0$. Then inequality (\ref{e4.2}) of Lemma~\ref{l4.2} yields
$$
\|u_n\|^2\leq C(\|u_n\|+\varepsilon_n^{1/2}\|u_n\|^{3/2}+\varepsilon_n\|u_n\|^2)\,.
$$
This implies the boundedness of $u_n$.

If $c=0$, then the boundedness of $u_n$ and inequality (\ref{e4.2}) imply that $\|u_n\|\to 0$. The converse implication is trivial.

\hfill$\Box$

\begin{prop}\label{p5.2}
Every Palais-Smale sequence for $J|_{\mathcal{N}}$ is a Palais-Smale sequence for $J$.
\end{prop}

{\em Proof\/}. Let $u_n\in\mathcal{N}$ be a Palais-Smale sequence for $J|_{\mathcal{N}}$. Inequality (\ref{e4.2}) of Lemma~\ref{l4.2} implies immediately that the sequence $u_n$ is bounded. Let $g_n=\nabla J(u_n)$ and $g_n^\tau$ be the tangent component of $g_n$, {\em i.e.}, the orthogonal projection of $g_n$ onto the tangent space   to $\mathcal{N}$ at $u_n$. Then, by assumption, $g_n^\tau\to 0$. We have to show that, actually, $g_n\to 0$.

Since the sequence $u_n$ is bounded, then, by Lemma~\ref{l4.3a}$(b)$, $\|P_n\|\leq C$ for some $C>0$, where $P_n$  is the projector onto ${\bar E}_{u_n}$ along the tangent space $T_{u_n}$. The  adjoint operator, $P_n$, is the projector onto the orthogonal complement to ${\bar E}_{u_n}$ along the normal subspace to $\mathcal{N}$ at $u_n$, and $\|P^*_n\|\leq C$ for some $C>0$ independent of $n$.

Now notice that, by the definition of $\mathcal{N}$, $g_n$ is orthogonal to the subspace ${\bar E}_{u_n}$. Therefore, $g_n=P^*_ng_n^\tau$ and, hence, $\|g_n\|\leq C\|g_n^\tau\|\to 0$. This completes the proof.

\hfill$\Box$

\vspace{1ex}

\begin{prop}\label{p5.3}
If $u_n$ is a Palais-Smale sequence for $J$ at a level $c>0$, then there exists a Palais-Smale sequence ${\tilde
u}_n\in \mathcal{N}$ at the same level such that $\|u_n - {\tilde u}_n\|\to 0$, and $c\geq m$.
\end{prop}

{\em Proof\/}. By Proposition~\ref{p5.1}, the sequence $u_n$ is bounded.  Hence,
$$
J(u_n)-\frac{1}{2}(J'(u_n),u_n)=\frac{1}{2}\int_\mathbb{R}f(x, u_n)u_n dx -\int_\mathbb{R}F(x,u_n) dx \to c\,.
$$
Since $F(x,u)\geq 0$, we obtain that
$$
\int_\mathbb{R}f(x, u_n)u_n dx\geq c
$$
for all $n$ large enough. Furthermore, Palais-Smale property also implies that  $G(u_n)\to 0$. By Lemma~\ref{l4.4},
$$
{\tilde u_n}= G^{-1}_{u_n}(0)\in \mathcal{N}
$$
is well-defined for all $n$ large enough, and $\|u_n - {\tilde u}_n\|\to 0$. Obviously, ${\tilde u}_n$ is a Palais-Smale
 sequence at the level $c$. Since $J({\tilde u}_n)\geq m$, the last statement of proposition follows immediately.

\hfill$\Box$

Combining Propositions~\ref{p4.1} and \ref{p5.3}, we obtain

\begin{cor}\label{c5.1}
If $u_n$ is a Palais-Smale sequence for $J$ at a positive level, then $\liminf \|u_n\|_\infty\geq\varepsilon_0>0$,
where $\varepsilon_0$ is the constant from Proposition~\ref{p4.1}.
\end{cor}

The following proposition is one of our key ingredients.

\begin{prop}\label{p5.4}
Given $\varepsilon>0$, there exists a Palais-Smale sequence $u_n$ for the functional $J|_\mathcal{N}$ (hence, for $J$)
at some level $c\in [m, m+\varepsilon]$ such that
\begin{equation}\label{e5.1}
\lim_{n\to\infty}\|u_{n+1}-u_n\|=0\,
\end{equation}
{\em i.e.}, a $(\overline{PS})$ sequence.
\end{prop}

{\em Proof\/}.
On the manifold $\mathcal{N}$ we consider the following initial-value problem
\begin{equation}\label{e5.2}
\frac{d\zeta}{dt}=-\nabla_\tau J(\zeta)\,,\quad \zeta(0)=u_0\in\mathcal{N}\,.
\end{equation}
The right-hand side of the differential equation in (\ref{e5.2}) is locally bounded and Lipschitz continuous.  Since 
$\mathcal{N}$ is a $C^{1,1}$ manifold, we see that problem (\ref{e5.2}) has a unique local solution 
$\zeta(t)\in\mathcal{N}$ for any $u_0\in\mathcal{N}$. Indeed, the problem reduces to an initial-value problem on a ball 
in $T_{u_0}$ centered at $0$, with Lipschitz continuous right-hand side. Moreover, if $\|u_0\|\leq R$ for some $R>0$, 
then, by Proposition~\ref{p4.2}, both the radius of the ball on which the reduced problem is defined and the Lipschitz 
constant of the right-hand side depend only on $R$, not on $u_0$. This implies that the local solution is defined on a 
time interval whose length is bounded below by a positive constant that depends only on $R$.

If $\zeta(t)$ is a solution of (\ref{e2.5}), then $J(\zeta(t))$ is a non-increasing function of $t$. Hence,
$J(\zeta(t))\leq J(u_0)$ for all positive $t$ in the domain of the solution. Therefore, by inequality (\ref{e4.2}) of
Lemma~\ref{l4.2}, $\|\zeta(t)\|\leq R$ on the domain of $\zeta$, where $R>0$ depends only on $J(u_0)$. This implies that
the solution is defined for all $t>0$.

Now we choose any $u_0\in\mathcal{N}$ such that $J(u_0)\leq m+\varepsilon$. Then
$$
J(\zeta(t))\to c\in [m, m+\varepsilon]
$$
as $t\to\infty$. Let
$$
\varphi(t)=-\int_0^t \|\nabla_\tau J(\zeta(s))\|^2 ds\,.
$$
Then
$$
\varphi(t)=J(\zeta(t))-J(u_0)
$$
and
$$
\lim_{t\to\infty}\varphi(t)=\inf_{t>0}\varphi(t)\geq -\varepsilon\,.
$$
Let $s_n\to\infty$ be a sequence such that $|s_n-s_{n-1}|\to 0$. Since $s_n$ is a minimizing sequence for $\varphi$,
Ekeland's variational principle implies the existence of a sequence $t_n$ such that $\varphi(t_n)\to
\inf_{t>0}\varphi(t)$, $\varphi'(t_n)\to 0$ and $t_n-s_n\to 0$. Setting $u_n=\zeta(t_n)$, we obtain a Palais-Smale
sequence for $J|_\mathcal{N}$. (Alternatively, at this point one can use an elementary argument from Real Analysis
instead of Ekeland's principle). Finally, since $\zeta(t)$ is bounded, $\nabla_\tau(\zeta(t))$ is bounded as well.
Hence, by the mean value theorem and the equation for $\zeta$,
$$
\|u_n-u_{n-1}\|\leq C |t_n-t_{n-1}| \to \infty\,.
$$
This completes the proof.

\hfill$\Box$

In what follows we consider Palais-Smale sequences at levels close to $m$. The next result shows that the
structure of such sequences is much simpler than in general case.

\begin{lemma}\label{l5.1}
Let $u_n$ be a Palais-Smale sequence for $J$ at a level $c\in [m, 2m)$. Suppose that $u_n\to u_0$ weakly in
$H^1(\mathbb{R})$.

$(a)$ If $u_0\neq 0$, then $u_n\to u_0$ strongly in $H^1(\mathbb{R})$, $u_0$ is a critical point of $J$, and $J(u_0)=c$.

$(b)$ If $u_0=0$, then there exist a sequence $x_n\in\mathbb{R}$, with $\lim |x_n|=\infty$, and a nontrivial critical
point $v_h$ of $J_h$  for some $h\in\mathcal{H}$, with $J_h(v_h)=c$, such that along a subsequence $T_{x_n}u_n\to v_h$
and
$u_n-T_{-x_n}v_h\to 0$ strongly in $H^1(\mathbb{R})$.
\end{lemma}

{\em Proof\/}. $(a)$ Since $J'$ is weakly continuous, $J'(u_0)=0$. By Proposition~\ref{p3a.1}, $u_n-u_0$ is a
Palais-Smale sequence at level $c-J(u_0)$. If $c-J(u_0)>0$, then, by Proposition~\ref{p5.3},  $c-J(u_0)\geq m$ which is
impossible because $c<2m$ while $J(u_0)\geq m$. Thus, $J(u_0)=c$ and $u_n-u_0$ is a Palais-Smale sequence at level
zero. By Proposition~\ref{p5.1}, $u_n-u_0\to 0$ strongly in $H^1(\mathbb{R})$.

$(b)$ Let $x_n\in\mathbb{R}$ be any point of global maximum for the function $|u_n|$ (obviously, such points exist), and
let $v_n=T_{x_n}u_n$. Since $u_n\to 0$ weakly, Corollary~\ref{c5.1} implies that $|x_n|\to \infty$.
Furthermore, zero is not a weak limit point of the sequence $v_n$. Since $v_n$ is a bounded sequence, then, along a
subsequence, $v_n\to v_h\neq 0$ weakly in $H^1(\mathbb{R})$. Passing to a subsequence one more time, we also obtain
limit potential $V_h$ and nonlinearity $f_h$, and, hence, the limit functional $J_h$.
By Proposition~\ref{p3.1} and weak continuity of $J'_h$, we obtain easily that $J'_h(v_h)=0$.

By Propositions~\ref{p3.1} and \ref{p3a.1},
\begin{equation*}
\begin{split}
J(u_n-T_{-x_n}v_h)-J(u_n)+J_h(v_h)=\\
= (J_{x_n}(v_n-v_1)-J_{x_n}(v_n)+J_{x_n}(v_h))-\\
- (J_{x_n}(v_h)-J_h(v_h))\to 0\,.
\end{split}
\end{equation*}
Hence, $J(u_n-T_{-x_n}v_h)\to c-J_h(v_h)$. Similarly, making use of second parts of Propositions~\ref{p3.1} and
\ref{p3a.1} we see that $J'_h(u_n-T_{-x_n}v_h)\to 0$ in $H^1(\mathbb{R})$ and, hence, $u_n-T_{-x_n}v_h$ is a
Palais-Smale sequence for $J$ at the level $c-J_h(v_h)$. As in the proof of first part of Proposition, we see that
$J_h(v_h)=c$ and, by Proposition~\ref{p5.1},
$$
\|u_n-T_{-x_n}v_h\|=\|v_n-v_h\|\to 0\,.
$$
This completes the proof.

\hfill $\Box$

\section{Proof of Main Result}\label{s6}

Theorem~\ref{t3.1} is an immediate consequence of the following two propositions.

\begin{prop}\label{p6.1}
If $u\in H^1(\mathbb{R})$ is a nontrivial solution of equation (\ref{e3.1}), then $u'$ is a continuous function and
$$
0<|u(x)|+|u'(x)|\leq \alpha\exp(-\beta|x|)
$$
for some positive constants $\alpha$ and $\beta$.
\end{prop}

{\em Proof\/}. Let $u\in H^1(\mathbb{R})$ be a nonzero solution. Set $V_1(x)=f(x,u(x))/u(x)$ for all $x$ such that
$u(x)\neq 0$ and $V_1(x)=0$ otherwise.  Then the function $u$ is an
$L^2$-eigenfunction of the operator $L-V_1(x)$ with the eigenvalue zero. It is easily seen that $V_1\in
L^\infty(\mathbb{R})$ and $V_1(x)$ vanishes at
infinity in the sense that $\mathrm{ess}\sup_{|x|\geq R}\to 0$ as $r\to\infty$. Hence, $V_1(x)$ is a relatively compact
perturbation of the operator $L$, and outside of $\sigma(L)$ the perturbed operator may have only isolated eigenvalues
of finite multiplicity. Now the result follows immediately from well-known properties of eigenfunctions of
Schr\"odinger operators (see, {\em e.g.}, \cite{si}).

\hfill$\Box$

\begin{prop}\label{p6.2}
For every $\varepsilon>0$ there exists a critical point of the functional $J$ with critical value $c\in [m, m+\varepsilon]$.
\end{prop}

{\em Proof\/}. Without loss of generality, we suppose that $\varepsilon<m$. Let $u_n$ be the  Palais-Smale sequence
from Proposition~\ref{p5.4}. We consider two cases.

\vspace{1ex}

{\bf Case 1}. The sequence $u_n$ has a non-zero weak limit point $u_0$. Then, by Lemma~\ref{l5.1}$(a)$, $u_0$ is a
critical point of $J$ at the level $c$, and we obtain the required.

\vspace{1ex}

{\bf Case 2}. The sequence $u_n$ converges to zero weakly in $H^1(\mathbb{R})$. For $u\in H^1(\mathbb{R})$ we set
$$
r(x; u)=\int_x^\infty[(u')^2(z)+u^2(z)]\,dz\,.
$$
The function $r(x; u)$ is continuous, non-increasing, and $r(x; u)\to 0$ as $x\to\infty$.  By
Proposition~\ref{p4.1},
there exists  $x_n\in \mathbb{R}$ such that $r(x_n, u_n)=\delta_0$, where $\delta_0=\varepsilon_0^2/2$. Note that $x_n$
is not necessarily unique.

We claim that $x_n-x_{n-1}\to 0$ and $|x_n|\to \infty$.
Consider any subsequence $x_{n'}$ of $x_n$. By Lemma~\ref{l5.1}$(b)$, there
exists a subsequence $u_{n''}$, a sequence $y_{n''}\in\mathbb{R}$, with $\lim |y_{n''}|=\infty$, and a nontrivial
critical point $v_h$ of $J_h$, for some $h\in\mathcal{H}$, such that $T_{y_{n''}}u_{n''}\to v_h$ strongly in
$H^1(\mathbb{R})$. Hence, $T_{y_{n''}}u_{n''-1}\to v_h$ strongly in
$H^1(\mathbb{R})$. This implies that
$$
r(x; T_{y_{n''}}u_{n''})\to r(x; v_h)
$$
and
$$
r(x; T_{y_{n''}}u_{n''-1})\to r(x; v_h)
$$
in $L^\infty(\mathbb{R})$. By Proposition~\ref{p6.1}, the function $r(x; v_h)$ is strictly decreasing, and there
exists a unique $x_h\in \mathbb{R}$ such that $r(x_h; v_h)=\delta_0$. Now it is easily seen that
$$
\lim x_{n''}-y_{n''}=\lim x_{n''-1}-y_{n''}= x_h\,.
$$
This implies the claim immediately.

Setting $v_n=T_{x_n}u_n$, we show that $0$ is not a weak limit point of the sequence $v_n$. Indeed, assume the
contrary. Since $u_n$ is a Palais-Smail sequence for $J$, then, along a subsequence, $v_n$ is a Palais-Smale sequence
for some functional in the envelope of $J$, and $v_n\to 0$ weakly in $H^1(\mathbb{R})$. By Lemma~\ref{l5.1}$(b)$,
passing to a further subsequence, there exists a sequence $y_n$,  yet another functional $J_h$
and its nontrivial critical point $v_h$ such that $|y_n|\to \infty$, $T_{y_n}v_n\to v_h$ and $T_{-y_n}v_h-v_n\to 0$
strongly in $H^1(\mathbb{R})$. This implies that
$$
r(0; v_n)-r(0;T_{-y_n}v_h)=\delta_0-r(0;T_{-y_n}v_h)\to 0\,.
$$
On the other hand, since $|y_n|\to \infty$, we see that, along a further subsequence, either
$$
r(0;T_{-y_n}v_h)\to 0\,,
$$
or
$$
r(0;T_{-y_n}v_h)\to \|v_h\|^2\geq 2\delta_0\,,
$$
a contradiction.

Now suppose for definiteness that the sequence $x_n$ is unbounded
above (the other case being similar). As it is well-known (see, {\em e.g.}, \cite{le-zhi, pa-ap}), there exists a
returning sequnce $z_k\to \infty$ for almost periodic functions $V(x)$ and $f_u(x,u)$ in the sense that
$$
T_{z_k}V\to V
$$
in $BS(\mathbb{R})$ and
$$
\|T_{z_k}f_u(\cdot, u)- f_u(\cdot, u)\|_R\to 0\quad \forall R>0\,.
$$
Since $x_n-x_{n-1}\to 0$, there exists a subsequence $x_{n_k}$ such that $z_k-x_{n_k}\to 0$. Along a subsequence,
$v_{n_k}\to v\neq 0$ weakly in $H^1(\mathbb{R})$. We shall show that $v$ is a solution of the problem. By
Proposition~\ref{p3.1},
\begin{equation*}\begin{split}
&\|J'(v_{n_k})-J'_{x_{n_k}}(v_{n_k})\|_*\leq \\
\leq&\|J'(v_{n_k})-J'_{z_{k}}(v_{n_k})\|_*+\|J'_{z_k}(v_{n_k})-J'_{x_{n_k}}(v_{n_k})\|_*\to 0\,.
\end{split}
\end{equation*}
Then the weak continuity of $J'$ and the fact that $u_n$ is a Palais-Smale sequence for $J$ imply that
\begin{equation*}\begin{split}
& (J'(v),\varphi)=\lim (J'(v_{n_k}),\varphi)  =\\
=& \lim (J'_{x_{n_k}}(v_{n_k}),\varphi)=\lim (J'(u_{n_k}),T_{-x_{n_k}}\varphi)=0\,.
\end{split}
\end{equation*}
This completes the proof.

\hfill$\Box$

\section{An Application}\label{s7}

Suppose that a dielectric medium occupies the whole space $\mathbb{R}^3$, and its material characteristics depend on
the $x$-variable only. Notice that in such a medium the magnetic permeability is equal to $1$ and, hence, the magnetic
induction is equal to  the magnetic field. Considering electromagnetic fields that depend on the time $t$ and the
$x$-variable only, we concentrate on a special class of such fields, the so-called $TE$-modes. In a $TE$-mode the
electric and magnetic components are of the form $(E, 0, 0)$ and $(0, H_1, H_2)$, respectively. Then the Maxwell
equations reduce to the following  equation
\begin{equation*}%\label{e7.1}
-\frac{\partial^2 E}{\partial x^2}=\frac{\partial^2\mathcal{F}(E)}{\partial t^2}
\end{equation*}
for the electric field only, where $D=\mathcal{F}(E)$ is the constitutive relation between the displacement and the
electric field.

In the so-called
Akhmediev-Kerr model this constitutive relation is of the form
$$
D=(\varepsilon(x)+g(x)\langle E^2 \rangle)E\,,
$$
where $\varepsilon(x)$ is the dielectric function of the medium, $g(x)$ represents nonlinear susceptibility, and
$\langle\cdot\rangle$ stands for the time average. In general, the nonlinear susceptibility may attain values of any
sign. If $g(x)>0$, the medium is self-focusing, while $g(x)<0$ means that the medium is defocusing.
However, in the following we assume that $g(x)$ does not change sign. More precicely, we suppose that the
functions $\varepsilon(x)$ and $g(x)$ are measurable, bounded, and Stepanov almost periodic, $g(x)$ does not change
sign, and both $\varepsilon(x)$ and $|g(x)|$ are bounded below by positive constants. Thus, we are dealing with a
one-dimensional almost periodic photonic crystal which is eather totally self-focusing, or totally defocusing,
depending on the sign of $g$.

A gap soliton is represented by a time-harmonic wave
$$
E=u(x)\cos(\omega t+\varphi_0)\,,
$$
where $\omega$ is a prohibited frequency and the wave profile $u(x)$ is a well-localized function vanishing at
infinity. After this Ansatz, we obtain the following equation for the profile function
$$
-\frac{d^2u}{dx^2}-\omega^2\varepsilon(x) u= g(x)u^3\,.
$$
Notice that the frequency $\omega$ is prohibited if and only if $0$ is not in the spectrum of the Schr\"odinger operator
$$
L=-\frac{d^2}{dx^2}-\omega^2\varepsilon(x)\,.
$$
Furthermore, since $\varepsilon(x)>0$, the operator $L$ is not positive definite, and the negative part of its spectrum
 is non-empty.

Thus, all the assumptions of Theorem~\ref{t3.1} are satisfied, and we obtain that there exists an exponentially
decaying wave profile $u\neq 0$. This shows that, in the framework of Akhmediev-Kerr model, one-dimensional almost
periodic photonic crystals possess gap solitons for all prohibited frequences.

\end{document}